\documentclass[12pt,dvipdfmx]{amsart}
\usepackage{a4}
\usepackage{latexsym}
\usepackage{amsmath}
\usepackage{amssymb}
\usepackage{amscd}
\usepackage{color}
\usepackage{amsthm}
\usepackage{hyperref}

\theoremstyle{plain}
    \newtheorem{theorem}                    {Theorem}[section]  
    \newtheorem{lemma}      [theorem]       {Lemma}
    \newtheorem{corollary}  [theorem]       {Corollary}
    \newtheorem{proposition}[theorem]       {Proposition}
    
    \newtheorem{conjecture} [theorem]       {Conjecture}

\newtheorem{example}[theorem]{Example}

\newtheorem{remark}[theorem]{Remark}

\newcommand{\Hom}{\operatorname{Hom}}
\newcommand{\End}{\operatorname{End}}

\newcommand{\Pic}{\operatorname{Pic}}

\newcommand{\disc}{\operatorname{disc}}
\newcommand{\tor}{{\operatorname{tor}}}

\newcommand{\Tor}{{\operatorname{Tor}}}

\newcommand{\im}{\operatorname{im}}

\newcommand{\coker}{\operatorname{coker}}

\newcommand{\Br}{\operatorname{Br}}

\renewcommand{\lim}{\operatorname{lim}}

\newcommand{\rk}{\operatorname{rank}}

\newcommand{\NS}{\operatorname{NS}}
\newcommand{\f}{{\mathcal F}}
\newcommand{\N}{{\mathbb N}}
\newcommand{\Z}{{{\mathbb Z}}}

\newcommand{\Q}{{{\mathbb Q}}}
\newcommand{\R}{{{\mathbb R}}}
\newcommand{\F}{{{\mathbb F}}}

\newcommand{\G}{{{\mathbb G}}}

\newcommand{\et}{{\text{\rm \'et}}}

\newcommand{\tr}{\operatorname{tr}}

\renewcommand{\O}{{\mathcal O}}

\newcommand{\proofend}{\hfill$\square$\\ \smallskip}

\date{\today}
\title{Brauer and N\'eron-Severi groups of surfaces over finite fields}
\author{Thomas H. Geisser}
\address{
	Department of Mathematics, Rikkyo University, Ikebukuro, Tokyo, Japan
}
\email{geisser@rikkyo.ac.jp}
\thanks{Supported by JSPS Grant-in-Aid (C) 18K03258}

\begin{document}

\begin{abstract}
We give a version of the Artin-Tate formula for surfaces over a finite 
field not assuming Tate's conjecture. It gives an equality between terms related 
to the Brauer group on the one hand and terms related to the N\'eron-Severi 
group on the other hand. We give estimates on the terms appearing in the formula
and use this to gives sharp estimates on the size of the Brauer group of abelian
surfaces depending on the $p$-rank. 
\end{abstract}

\maketitle

\section{Introduction}

Let $\F_q$ be a finite field with Galois group $G$ and Weil group $\Gamma$,
and $X$ be a smooth, proper, and geometrically connected surface over 
$\F_q$ for which Tate's conjecture holds,
or equivalently, such that the Brauer group $\Br(X)$ is finite.
Let $\bar X$ be the base extension to the algebraic closure
and let $P_2(X,t)=\det(1-Ft\mid H^2_\et(\bar X,\Q_l))=\prod_i(1-\alpha_it)$ 
be the second $L$-polynomial of $X$. 
Artin-Tate \cite{artintate} and Milne \cite{milne} showed that
$$ \frac{P_2(X,q^{-s})}{(1-q^{1-s})^{\rho(X)}}\sim 
\pm \frac{|\Br(X)\cdot |\det\NS(X)|}{q^{\alpha(X)} |\NS(X)_\tor|^2}
\qquad \text{as}\; s\to 1,$$
where $\alpha(X)=\chi(X,\O_X)-1+\dim \text{PicVar}_X$, 
$\NS(X)$ is the image of $\Pic(X)$ in $\NS(\bar X)$, 
the Picard number $\rho(X)$ is the rank of $\NS(X)$, and $\det\NS(X)$ is the 
determinant of the intersection pairing on $\NS(X)/\tor$. 
It is known that the order $|\Br(X)|$ of the Brauer group is a square
\cite{ichbrauer,lgr}. 

In this note, we prove a version of this formula which does not assume
the finiteness of the Brauer group and does not involve $P_2(X,t)$,
and give estimates on the terms in the formula occurring. There are not many 
non-trvial examples of surfaces over finite fields for which the size of the Brauer group is known.
We calculate Brauer groups of abelian varieties for arbitrary $\F_q$, and use
this to show that our estimates are sharp.

Consider the map 
$$\alpha_X:\Br(X)\to \Br(\bar X)^\Gamma.$$ 
Based on work of Colliot-Th\'el\`ene and Skorobogatov \cite{cts}, 
Yuan showed that the kernel and cokernel of $\alpha_X$ are finite 
\cite[Thm.\ 1.3(4)]{yuan}. Let 
$$\beta_X : \NS(\bar X)^\Gamma \to\NS(\bar X)_\Gamma$$
be the map from invariant to coinvariants induced by the identity, and let 
$\omega_X:=z(\beta_X)^{-1}$, where we set $z(f)=\frac{|\ker \beta|}{|\coker \beta|}$ 
for any map $f$ of abelian groups with finite kernel and cokernel.
Finally, let $H^2_W(X,\G_m )/\tor\stackrel{\bar d}{\to}  \Hom(H^1_W(X,\G_m ),\Z)$
be the map induced by the cup-product of Weil-\'etale cohomology groups. 
Then $|\coker \bar d|$ is a divisor
of $|\det \NS(X)|$, and it is equal to one if the Brauer group is finite.
In Section 2 we show the following equality between terms related to the Brauer
group and terms related to the N\'eron-Severi group:  

\begin{theorem}[\ref{mainth}] 
The map $\beta_X$ has finite kernel and cokernel, and $\omega_X\in \N$.
We have  
$$z(\alpha_X)=\frac{\omega_X\cdot |\NS(X)_\tor|\cdot |\coker \bar d|}{|\det \NS(X)|}.$$
\end{theorem}

\begin{corollary}[\ref{maincor}]
If $\Br(X)$ is finite, then 
$$ \frac{|H^1(G,\NS(\bar X))|}{|\coker\alpha_X|}=
\frac{|\Br(X)|}{|\Br(\bar X)^\Gamma|}
=\omega_X\frac{|\NS(X)_\tor|}{|\det \NS(X)|.}$$
\end{corollary}

If the Galois group acts trivially on $\NS(\bar X)$, 
then %$\omega=1$ and  
this simplifies to
$$ |\det \NS(X)|= |\coker\alpha_X|.$$

The theorem raises the question about the possible values of $\omega_X\in \N$
and the size of $\Br(\bar X)^\Gamma$.
Let $\rho(\bar X)=\rk \NS(\bar X)$ and $\rho(X)=\rk \NS(X)$ be the 
geometric and arithmetic Picard numbers, respectively, and $\tau_X$ its difference.
In Section 3 we prove:

\begin{proposition}[\ref{omegaX}] For any prime number $p$, 
the power of $p$ appearing in $\omega_X$ is at most 
$\frac{\tau_X}{p-1}$. In particular, the prime divisors of $\omega_X$ and are 
not larger than $\tau_X+1$, and if $\tau_X+1$ is a prime divisor of $\omega_X$, 
then $\omega_X=\tau_X+1$.
\end{proposition}

For example, the possible values are $\omega_X=2$ for $\tau_X=1$, 
$\omega_X\in\{1,2,3,4\}$ for $\tau_X=2$, $\omega_X\in \{2,4,6,8\}$ for $\tau_X=3$,
and $\omega_X\in \{1,2,3,4,5,6,8,9,12,16\}$ for $\tau_X=4$.

\smallskip

To describe $\Br(\bar X)^\Gamma$, we 
use results and methods of Milne \cite{milne}: Let $L_{tr}(X,t)$ be the 
transcendental part of $P_2(X,t)$, i.e., the product $\prod(1-\alpha_i t)$
over those $\alpha_i$ which are not $q$ times a root of unity.
Let $2g$ be its degree. In Section 4 we prove:

\begin{theorem}[\ref{mainbr}]
If $\Br X$ is finite, then
$$ |\Br(\bar X)^\Gamma|= |\NS(X)_\tor|\cdot 
L_{tr}(X,q^{-1})\cdot q^{\alpha(X)}.$$
\end{theorem}

We estimate the size of the $L$-value in Section 5.

\begin{proposition}[\ref{g=2}]
We have $$|L_{tr}(X,q^{-1})|\leq (4-\frac{1}{q})^g.$$
%In particular, if $\Br^0(\bar X)$ is the divisible part of the Brauer group, 
%then
%$$|\Br^0(\bar X)^\Gamma|\leq (4q-1)^g.$$
\end{proposition}

We deduce this from the following fact about Weil-polynomials, which in turn follows
from classical results on estimates of totally positive algebraic integers:

\begin{proposition}[\ref{c2}]
Assume $f(x)=x^{2g}+a_1x^{2g-1}+\cdots +a_1q^{2g-1}x+q^{2g}$ is a 
Weil-$q^2$-polynomial.
If $a_1>(2q-1)g$, then $f(-q)=0$.
\end{proposition}

The estimate is sharp because $(x^2+(2q-1)x+q^2)^g$ has $a_1=(2q-1)g$ and 
does not vanish at $-q$. 

We use the above result to give a complete picture in case of abelian surfaces,
building on calculations of Suwa-Yui \cite{suwayui}.
Among others we show: 

\begin{proposition}[\ref{ss}, \ref{largest}, \ref{largest2}]
a) For any abelian surface over $\F_q$, we have 
$$ |\Br(A)|\cdot|\det \NS(A)|\leq 16 q,$$
and for every $q$ there is an abelian surface where equality holds.
There is an abelian surface with $|\Br(A)|=16q$ if and only if $q$ is a square.

b) If $A$ is not supersingular, then
$$ |\Br(A)|\cdot |\det \NS(A)|\leq 16q-4$$
and for every $q$ there is an $A$ such that equality holds.

c) If $A$ has $p$-rank one, then 
$$|\Br(A)|\cdot |\det \NS(A)|\leq (4\sqrt q-1)^2$$
There is an abelian surface of 
$p$-rank one with $|\Br(A)|=(4\sqrt q-1)^2$ if and only if $q$ is a square.
\end{proposition} 

This gives a clear picture for the order of the Brauer group if $q$ a square, but 
the situation for 
$q$ not a square is more complicated. We have the following result which
illustrates that the order of the Brauer group does not necessarily increase
with the prime field:

\begin{proposition}[\ref{ssp}]
The abelian surface $A$ over $\F_3$ with Weil-polynomial 
$x^4-3x^2+9$ %a product of two non-isogenous curves, 
has $|\Br(A)|=36$, 
and this is the largest order of a Brauer group for a supersingular abelian 
surface over a prime field. 
\end{proposition}

\medskip

\noindent{\bf Notation:} We let  $\Q/\Z'=\Q/\Z[\frac{1}{p}]$. 
For any map $f$ of abelian groups 
with finite kernel and cokernel we set $z(f)=\frac{|\ker f|}{|\coker f|}$. 
For an abelian group $A$ we let $A_\tor$ the torsion subgroup and 
$A/\tor$ the torsion free quotient of $A$. We frequently switch between Weil-$q^2$-polynomials
$x^{2g}+a_1x^{2g-1}+\cdots +a_1q^{2g-2}x+q^{2g}$ and the corresponding
$L$-polynomial $1+a_1t+\cdots +a_1q^{2g-2}t^{2g-1}+q^{2g}t^{2g}$.

\section{The main theorem}
Let $X$ be a smooth, proper, and geometrically connected surface over a finite 
field $\F_q$. %for which the Tate conjecture holds 
Let $G$ be the Galois group and $\Gamma$ be the 
Weil group of $\F_q$ with Frobenius $F$ as the generator. For a G-module $A$,
we note  that $A^\Gamma=A^G$.  
The intersection pairing induces an exact sequence 
$$0\to \NS(X)_\tor\to \NS(X)\to \Hom(\NS(X),\Z)\to D\to 0$$
such that $|D|=|\det \NS(X)|$. Consider the commutative diagram
$$\begin{CD}
\Pic(X)@>>>\NS(X)\\
@VVV@VVV\\ 
\Pic(\bar X)^G@>>> \NS(\bar X)^G.
\end{CD} $$
By  definition the upper horizontal map is surjective and the right vertical map
is injective. The lower horizontal map is surjective by Lang's theorem  
$H^i(G,\Pic^0(\bar X))=H^i(\Gamma,\Pic^0(\bar X))=0$ for $i\geq 1$.  
%The left vertical map is an isomorphism because 
%$\Pic^0(X)\cong \Pic^0_X(\F_q)\cong \Pic^0_X(\bar \F_q)^G\cong \Pic^0(\bar X)^G$. 
The left vertical map is an isomorphism by
the Hochschild-Serre spectral sequence for $\G_m$ because 
the kernel $H^1(G, H^0(\bar X,\G_m))$ vanishes by Hilbert 90, and the cokernel is
$H^2(G, H^0(\bar X,\G_m))\cong \Br(\F_q)=0$. 
Hence $\NS(X)\cong \NS(\bar X)^G$. 

%The action of $G$ and $\Gamma$ on $\NS(\bar X)$ factor through a finite quotient.
Lang's theorem implies that $H^i(G,\Pic^0(\bar X))\cong H^i(G,\NS(\bar X))$ for $i\geq 1$,
hence the low terms of the Hochschild-Serre spectral sequence 
$H^s(G,H^t_\et(\bar X,\G_m))\Rightarrow H^{s+t}_\et(X,\G_m)$ for the cohomology of $\G_m$ 
give an exact sequence
$$ 0\to H^1(G,\NS(\bar X))\to \Br(X)\stackrel{\alpha_X}{\longrightarrow} 
\Br(\bar X)^\Gamma \to H^2(G,\NS(\bar X)) \to H^3(X,\G_m),$$ 
and the cokernel of $\alpha$ is finite \cite{yuan}, so that 
\begin{equation}\label{alphavalue}
z(\alpha_X)=\frac{|H^1(G,\NS(\bar X))|}{|\coker\alpha_X|},
%= \frac{|\Br(X)|}{|\Br(\bar X)^\Gamma|},
\end{equation}
and this is equal to $\frac{|\Br(X)|}{|\Br(\bar X)^\Gamma|}$
if $\Br(X)$ is finite.

We are using Weil-\'etale cohomology of $\G_m$, 
see \cite{Geisser04} for a definition and basic properties. 
For any bounded complex $\f$ of \'etale sheaves 
we have a long exact sequence \cite[Thm. 7.1c]{Geisser04}
$$ \cdots\to H^i_\et(X,\f)\to H^i_W(X,\f) \to H^{i-1}_\et(X,\f)_\Q
\stackrel{\partial}{\to} H^{i+1}_\et(X,\f)\to \cdots$$
which gives $H^1_W(X,\G_m )\cong \Pic(X)$ and $H^2_W(X,\G_m )_\tor\cong\Br(X)$.
The map $\partial$ can be written by \cite[Thm.\ 5.3]{Geisser04} as the composition
$$ \begin{CD}
H^{i-1}_\et(X,\f)_\Q@>\partial>> H^{i+1}_\et(X,\f)\\
@VVV @AAA\\
H^{i-1}_\et(X,\f\otimes^L\Q/\Z)@>\epsilon>> H^{i}_\et(X,\f\otimes^L\Q/\Z),
\end{CD}$$
where $\epsilon$ is the cup product with a generator of $H^1_W(\F_q,\Z)\cong \Z$.
In the case $i=2$ and $\f=\G_m$ this implies that the image of $H^2_W(X,\G_m )/\tor$
in $\NS(X)_\Q=H^1_\et(X,\G_m)_\Q$, i.e., the kernel of $\partial$,
contains $\NS(X)/\tor=H^1_\et(X,\G_m)/\tor$, i.e., the kernel of the left vertical map,
as a sublattice of some index $I_X$.
Cup product induces the map in the middle of the following commutative diagram,
in which the left groups are lattices in the right groups: 
\begin{equation}\label{pairing}
\begin{CD}
\NS(X)/\tor @>>> H^2_W(X,\G_m )/\tor@>>> \NS(X)_\Q \\
@VVV @V\bar d VV @V d_\Q V\simeq V \\
\Hom(\NS(X),\Z)@<\simeq<<  \Hom(H^1_W(X,\G_m ),\Z)@>>> \Hom(\NS(X)_\Q,\Q) .
\end{CD}
\end{equation}
Since the intersection pairing is non-degenerate after tensoring with 
$\Q$, the right square shows that $\bar d$ is injective and the left square
shows that
\begin{equation}\label{IX}
|\coker \bar d|\cdot I_X=|\det \NS(X)|.
\end{equation}

\begin{lemma} The number $\omega_X:=z(\beta_X)^{-1}$ is an integer. 
\end{lemma}

\proof
The Galois group acts through a finite quotient on $\NS(\bar X)$,
so that $\beta_X$ has finite kernel and cokernel.
Set $N=\NS(\bar X)$. 
Applying the snake Lemma to the self-map $F-1$ of the short
exact sequence
$0\to N_\tor \to N  \to (N/\tor)\to 0$
we see that in the double complex the kernel of the right upper map is 
isomorphism to the kernel of the lower left map.
\begin{equation}\label{nstor}
\begin{CD} 
(N_\tor)^\Gamma @>>> N^\Gamma @>>> (N/\tor)^\Gamma\\
@V\beta^\tor_X VV @V\beta_X VV@V\beta'_XVV\\
(N_\tor)_\Gamma @>>> N_\Gamma @>>> (N/\tor)_\Gamma\\
\end{CD}
\end{equation}
Hence the double complex is quasi-isomorphic to zero. 
Since $z(\beta^\tor_X)=1$ we obtain that $z(\beta_X)=z(\beta'_X)$
and the result follow because $\beta'_X$ is injective. 
\proofend

%$$ H^i_W(X.\G_m )_\tor \times H^{4-i}_W(X,\G_m )_\tor \to \Q/\Z.$$

\begin{theorem}\label{mainth} 
We have  
%$$z(\alpha_X)=\omega_X\frac{|\NS(X)_\tor|}{|\det \NS(X)|}\cdot |\coker \bar d|.$$
$$z(\alpha_X)=\frac{\omega_X\cdot |\NS(X)_\tor|}{I_X}.$$
%If $G$ acts trivially on $\NS(\bar X)$, then  
%$$ |\det \NS(X)|= |\coker\alpha_X|\cdot |\coker \bar d|.$$
\end{theorem}

\proof
The short exact sequence \cite[(9)]{Geisser04}
$$ 0\to H^{i-1}_\et(\bar X,\f)_\Gamma\to H^i_W(X,\f)\to 
H^i_\et(\bar X,\f)^\Gamma\to 0$$
gives $H^1_W(X,\G_m )\cong \Pic(\bar X)^\Gamma$ and a short exact sequence
\begin{equation}\label{h3w}
0\to \NS(\bar X)_\Gamma\stackrel{j}{\to} H^2_W(X,\G_m )\to \Br(\bar X)^\Gamma\to 0
\end{equation}
since $\Pic^0(\bar X)_\Gamma=0$. 
Consider the commutative diagram in which the map $d$ is induced by $\bar d$: 
$$\begin{CD}
\NS(X) @>e>> \Hom(\NS(X),\Z) @<g^* <\sim<\Hom(H^1_W(X,\G_m ),\Z)\\
@| @. @AdAA\\
\NS(\bar X)^\Gamma @>\beta_X >> \NS(\bar X)_\Gamma @>j>> H^2_W(X,\G_m ).
\end{CD}$$
The map $g^*$ is an isomorphism because $H^1_W(X,\G_m )\cong \Pic X$
and $\Hom(A,\Z)\cong \Hom(A/\tor,\Z)$ for any abelian group. We obtain
$ z(\beta_X)z(dj)=z(e)=\frac{|\NS(X)_\tor|}{|\det \NS(X)|}$. 
Since $j$ is injective, we have an exact sequence
$$ 0\to \ker dj\to \ker d\to \coker j\to \coker dj\to \coker d\to 0.$$
From the diagram \eqref{pairing}
we see that $\ker d=H^2_W(X,\G_m )_\tor \cong\Br(X)$
and that the map $\Br(X)\cong \ker d\to \coker j\cong  \Br(\bar X)^\Gamma$ 
is the natural map $\alpha_X$. 
Thus, $z(\alpha_X)=z(dj)\cdot |\coker \bar d|$ and we obtain the formula 
with \eqref{IX}.
\proofend

\begin{corollary}\label{maincor}
If $\Br(X)$ is finite, then 
$$ \frac{|H^1(G,\NS(\bar X))|}{|\coker\alpha_X|}=
\frac{|\Br(X)|}{|\Br(\bar X)^\Gamma|}
=\omega_X\frac{|\NS(X)_\tor|}{|\det \NS(X)|.}$$
If moreover $G$ acts trivially on $\NS(\bar X)$, then  
$$ |\det \NS(X)|= |\coker\alpha_X|.$$
\end{corollary}

\proof 
The pairing $\bar d$ is perfect if $\Br(X)$ is finite,
see \cite[Remark 3.3]{tata} or the proof of \cite[Cor. 4.2]{curve}
because finiteness of $\Br(X)$ implies that  all groups $H^i_W(X,\G_m )$ are 
finitely generated. Thus it suffices to combine the theorem with 
\eqref{alphavalue}. Finally, 
if $\Gamma$ acts trivially on $\NS(\bar X)$, 
then $\omega_X=1$ and $H^1(G,\NS(\bar X))\cong \NS(X)_\tor$.
\proofend

\section{Estimating $\omega_X$}
In this section we give estimates for $\omega_X$ in terms of 
$\tau_X=\rk \NS(\bar X)-\rk \NS(X)$ and $n$,
the degree of the minimal extension such that the Galois group acts trivially on 
$\NS(\bar X)$.
For an integer $b$, let $\varphi(b)$ be the Euler $\varphi$ function,
and $\lambda(b)=p$ if $b=p^t$ is a prime power and $\lambda(b)=1$
otherwise, i.e., $\lambda(b)=\exp(\Lambda(b))$ where $\Lambda$ is the von 
Mangoldt function. 

\begin{theorem}
%Let $\tau :=\rho(\bar X)-\rho(X)$, and 
There are non-negative integers $a_d$ for $d|n$ such that
$$\omega_X=\prod_{d|n}\lambda(d)^{a_d},\qquad 
\sum_{1\not=d|n} a_d\varphi(d)=\tau_X,\qquad 
n=\mathrm{lcm}\{d\mid a_d>0\}.$$
\end{theorem}

To prove the theorem, we can replace $\NS(\bar X)$ by its maximal torsion 
free quotient $N$ 
by the argument in \eqref{nstor}, and we let $M$ be the quotient of $N$ by 
image of the saturated submodule 
$\NS(X)\cong \NS(\bar X)^\Gamma$; note that $M\cong \im(F-1)$, hence $M$ is free,
of rank $\tau_X$.

\begin{lemma} 
We have $\omega_X=|\det(F-1|_{M_\Q})|$.
\end{lemma}

\proof %We use the method of \cite[Lemma z.4]{artintate}:
Let $f'$ be the restriction of $F-1$ to $M$. Then  
\begin{align*}
\ker \beta&=\ker( F-1)\cap\, \im( F-1)=\ker f';\\
\coker \beta&=N/(\ker( F-1)+\im (F-1))\cong M/\ker f'\cong \coker f'.
\end{align*}
Since $M$ is finitely generated free, we can apply the elementary divisor theorem
and assume that $f'$ is represented by a diagonal matrix $A$. 
In this case $f'$ is injective
because $\beta$ is a rational isomorphism, and $|\coker f'|=|\det(A)|$. 
The determinant does not change if we extend the base  ring from $\Z$ to $\Q$.
\proofend

Since $F$ acts via the finite cyclic group $\Z/n\Z$ on $M_\Q$, 
%with $n$ the degree of the 
%minimal extension such that the Galois group acts trivially on $\NS(\bar X)$.
$M_\Q$ is a direct sum of $\Q$-simple representations of $\Z/n\Z$.
If we write $x^n-1=\prod_{d|n}\Phi_d(x)$ as the product of cyclotomic
polynomials, then those simple representations are isomorphic to 
$\Q[x]/(\Phi_d(x))$, for $d|n$. Here $F$ acts as multiplication by $x$, and
$\dim \Q[x]/(\Phi_d(x))=\varphi(d)$.
Since $M_\Q$ does not contain the trivial representation corresponding to $d=1$,
we obtain that for some non-negative integer $a_d$ we have
\begin{equation}\label{ds}
M_\Q\cong \bigoplus_{1\not=d|n}\Q[x]/(\Phi_d(x))^{a_d}, \qquad 
\sum a_d\varphi(d)=\tau_X,
\end{equation}
proving the middle formula.
We get the formula for $n$ by noting
that the order of $F$ on $\Q[x]/(\Phi_d(x))$ is $d$, hence 
$a_d>0$ implies that $d|n$, and $n$ is the smallest number with this property.

Finally, it is easy to see with the basis $\{1,x,\ldots, x^{\varphi(d)-1}\}$
that the characteristic polynomial of $x$ acting on 
$\Q[x]/(\Phi_d(x))$ is $\Phi_d(t)$.
Thus the determinant of $F-1$ is $\Phi_d(1)$, and it remains to show: 

\begin{lemma}
$\Phi_d(1)=p$ if $d=p^i$ is a prime power, and $\Phi_d(1)=1$ otherwise.
\end{lemma}

\proof
Dividing the formula 
$\prod_{d|n} \Phi_d(x)=x^n-1$ %\frac{x^n-1}{x-1}=x^{n-1}+\cdots +x+1$ 
by $x-1$ and evaluating at $1$ we obtain $\prod_{1\not=d|n} \Phi_d(1)=n$.
The Lemma then follows by induction on $n$.
\proofend

\begin{remark} It follows from the proof that $\omega_X$ only depends on 
$\NS(\bar X)\otimes\Q$. However, 
$H^1(G,\NS(\bar X))\cong \Tor(\NS(\bar X)_\Gamma)$ depends on the 
integral structure. For example, let $M=\Z\oplus \Z$ with action 
$F(1,0)=(1,0), F(0,1)=(0,-1)$ and $N=\Z\oplus \Z$ with action
$F(1,0)=(0,1), F(0,1)=(1,0)$. Then $M$ is a sublattice on $N$ of index $2$,
$M_\Gamma \cong \Z\oplus\Z/2\Z$ and $H^i(G,M)\cong\Z/2\Z$ for all $i>0$
whereas $N_\Gamma\cong \Z$ and $H^i(G,N)=0$ for all $i>0$. 
\end{remark}

The theorem gives restrictions on $\omega_X$ and $n$ in terms of $\tau_X$.

%\begin{corollary}
%Assume that $n$ is a power of $p$. Then $p-1|\tau$
%and $\omega=p^{\frac{\tau}{p-1}}$.
%\end{corollary}

\begin{corollary}\label{omegaX}
1) The prime divisors of $\omega_X$ divide $n$.

2) The power of $p$ appearing in $\omega_X$ is at most $\frac{\tau_X}{p-1}$. 
In particular,
the prime divisors of $\omega_X$ and are not larger than $\tau_X+1$, and if 
$\tau_X+1$ is a prime divisor of $\omega_X$, then $\omega_X=\tau_X+1$.
\end{corollary}

\proof
The first statement is clear. 
If $p|\omega_X$, then $M_\Q$ contains a subrepresentation 
$\Q[x]/(\Phi_{p^r}(x))$. 
But then $\tau_X\geq \varphi(p^r)=p^r-p^{r-1}\geq p-1$. Equality can
only happen if $M_\Q\cong \Q[x]/(\Phi_{p}(x))$.
\proofend

\begin{corollary}\label{exactomega}
1) If $\tau_X$ is odd, then the possible values for $(n,\omega_X)$ are 
$(lcm(n',2), 2\omega_X')$ for $(n',\omega_X')$ a possible value for $\tau_X-1$.

2) For small $\tau_X$, we have the following possibilities for $n$ and $\omega_X$:
$$\begin{array}{r|r|r|r}
\tau_X &d's\; in\; \eqref{ds}& n & \omega_X\\
\hline
1&(2)&2&2\\
\hline
2&(2,2)&2&4\\
&(3)&3&3\\
&(4)&4&2\\
&(6)&6&1\\
\hline
3&(2,2,2)&2&8\\
&(2,3)&6&6\\
&(2,4)&4&4\\
&(2,6)&6&2\\
\end{array}%\quad
\begin{array}{r|r|r|r}
\tau_X &d's\; in\; \eqref{ds}&n & \omega_X\\
\hline
4&(2,2,2,2)&2&16\\
&(2,2,3)&6&12\\
&(2,2,4)&4&8\\
&(2,2,6)&6&4\\
&(3,3)&3&9\\
&(3,4)&12&6\\
&(3,6)&6&3\\
&(4,4)&4&4\\
&(4,6)&12&2\\
&(6,6)&6&1\\
&(5)&5&5\\
&(8)&8&2\\
&(10)&10&1\\
&(12)&12&1\\
\end{array}\quad
\begin{array}{r|r|r|r}
\tau_X & d's\; in\; \eqref{ds} &n & \omega_X\\
\hline
5&(2,2,2,2,2)&2&32\\
&(2,2,2,3)&6&24\\
&(2,2,2,4)&4&16\\
&(2,2,2,6)&6&8\\
&(2,3,3)&6&18\\
&(2,3,4)&12&12\\
&(2,3,6)&6&6\\
&(2,4,4)&4&8\\
&(2,4,6)&12&4\\
&(2,6,6)&6&2\\
&(2,5)&10&10\\
&(2,8)&8&4\\
&(2,10)&10&2\\
&(2,12)&12&2\\
\end{array}
$$

3) For $\tau_X=6$, $\omega_X$ is contained in the following set:
$$ 1,2,3,4,5,6,7,8,9,10,12,15,16,18,20,24,27,36,48, 64.$$
\end{corollary}

\proof
1) This follows because the only integer $d>1$ with $\varphi(d)$ odd is $d=2$,
hence in order for $\tau_X$ to be odd, $M_\Q$ needs to contain the 
representation of rank one where $F$ acts as $-1$.

2) We have to find a sequence of $d_j$ %  $\Q(\zeta_{d_j})$ 
such that $\sum\varphi(d_j)=\tau_X$. If $\tau_X=1$, then $n=2$ and $\omega_X=2$.
If $\tau_X=2$, then we can have two one-dimension representations ($n=2$) and get
$\omega_X=4$, or we have one two-dimensional representation and $n=3,4,6$
so that $\varphi(n)=2$. In this case we obtain $\omega_X=3, 2, 1$, respectively.
The case $\tau_X=3$ follows from (1), and the case $\tau_X=4$ is done similarly, 
using that $\varphi(n)=4$ for $n=5,8,10,12$ with $\omega_X=5,2,1,1$ respectively.

3) For $\tau_X=6$ one uses that $\varphi(n)=6$ only for $n=7,9,14,18$.
\proofend

\section{The invariants of the Brauer group}
Let $\Br^0(\bar X)$ be the maximal divisible subgroup of $\Br(\bar X)$, with quotient
the inverse limit $\Br^c(\bar X)=\lim \Br(\bar X)/m$.
%$\Br(\bar X)\cong \Br^0(\bar X)\oplus\Br^c(\bar X)$ as Galois modules. 
Recall that the functor $T\mapsto \underline H^3(X_T,\Z/p^r\Z(1))$ 
is represented 
by a group scheme on the perfect site (i.e., the category of perfect affine 
$\F_p$-schemes equipped with the \'etale topology) \cite{milne}, 
with connected quasi-algebraic subgroup $U^3(p^r)$ and \'etale 
quotient $D^3(p^r)$. 
The inverse system $U^3(p^r)$ becomes constant, and we let $s$ be its dimension.

\begin{proposition}\label{brautor}
We have 
$$| \Br^c(\bar X)^\Gamma| \cong |\NS(X)_\tor|\cdot q^s.$$
\end{proposition}

One way of proving this is to show that  
$\Br^c(\bar X)^\Gamma\cong H^2(X,\hat \G_m )^\Gamma$
and using \cite[Lemma 5.2]{milne}, but we give a self-contained proof. 

\proof
Away from $p$, $\NS(\bar X)_\tor$ is dual to 
$\Br^c(\bar X)$ by Poincar\'e duality, hence 
$|\Br^c(\bar X)^\Gamma|=|(\NS(\bar X)_\tor)_\Gamma|=|\NS(\bar X)_\tor^\Gamma|
=|\NS(X)_\tor|$. At the prime $p$, we have a diagram
$$ \begin{CD}
@. 0@. 0@. 0\\ 
@. @VVV@VVV@VVV\\
0@>>> A_r @>>> U^3(p^r)(\bar \F_q)@>>> B_r@>>>0\\
@. @VVV@VVV@VVV\\
0@>>> \Br(\bar X)/p^r 
@>>> H^3(\bar X, \Z/p^r\Z(1))@>>> H^3(\bar X,\G_m)[p^r]@>>>0\\
@. @VVV@VVV@VVV\\
0@>>> C_r @>>> D^3(p^r)(\bar \F_q)@>>> E_r @>>>0\\
@. @VVV@VVV@VVV\\
@. 0@. 0@. 0\\ 
\end{CD}$$
%We note that the systems $U^3(p^r)$, hence $A_r$ and $B_r$, become 
%constant for large enough $r$ and are groups of finite exponent. 
Taking the limit we obtain: 
$$ \begin{CD}
@. 0@. 0@. 0\\ 
@. @VVV@VVV@VVV\\
0@>>> \widehat A @>>> U^3(p^\infty)(\bar \F_q)@>>> \widehat B@>>>0\\
@. @VVV@VVV@VVV\\
0@>>> \Br(\bar X)^{\wedge p} 
@>>> H^3(\bar X, \Z_p(1))@>>> T_pH^3(\bar X,\G_m) @>>>0\\
@. @VVV@VVV@VVV\\
0@>>> \widehat C
@>>>\lim D^3(p^\infty)(\bar \F_q)@>>> \widehat E@>>>0\\
@. @VVV@VVV@VVV\\
@. 0@. 0@. 0\\ 
\end{CD}$$
Since $U^3(p^\infty)$ is of finite exponent, so is $\widehat B$. As 
$T_pH^3(\bar X,\G_m)$ is torsion free, we obtain $\widehat B=0$, 
hence $\widehat A\cong U^3(p^\infty)(\bar \F_q)$
and $T_pH^3(\bar X,\G_m)\cong \widehat E$. Thus the left horizontal arrows
are isomorphisms on torsion subgroups. Furthermore 
the group $D^3(p^r)(\bar \F_q)$ is dual to 
$D^1(p^r)(\bar \F_q)\cong \Pic(\bar X)[p^r]$, hence 
$\lim D^3(p^r)(\bar \F_q)\cong\Hom(\Pic(\bar X)\{p\},\Q/\Z)$ has torsion group
$\NS(\bar X)\{p\}^*$ and we obtain a short exact sequence
$$ 0\to U^3(p^\infty)(\bar \F_q)\to  \Br^c(\bar X)\{p\}\to 
\NS(\bar X)\{p\}^*\to 0.$$
Taking Galois invariants is exact because unipotent groups have no higher
cohomology, and we get
$$ 0\to U^3(p^\infty)(\F_q) \to \Br^c(\bar X)\{p\}^\Gamma\to 
(\NS(\bar X)\{p\}^*)^\Gamma\to 0.$$
\proofend

Write $P_2(X,t)=\prod (1-\alpha_i t)\in \Z[t]$ and let 
$$L_{tr}(X,t)=\prod_{\alpha_i\not=\zeta q} (1-\alpha_i t)\in \Z[t]$$
be its transcendental part, i.e., the product over those roots which
are not $q$ times a root of unity.  Since the $\alpha_i$ come
in complex conjugate pairs, the degree $2g$ of $L_{tr}(X,t)$ is even. 
%Moreover, $|\alpha_i|= q$ implies that 
%$(1-\alpha_iq^{-1})(1-\bar\alpha_iq^{-1})=2-\frac{\alpha_i+\bar \alpha_i}{q}$ 
Since $L_{tr}(X,t)$ has no real zero, $L_{tr}(X,q^{-1})$ is a positive real
number.

Let $T_l\Br(\bar X)$  be the Tate module and 
$V_l\Br(\bar X)=T_l\Br(\bar X)\otimes_{\Z_l}\Q_l$. 
From the exact sequence
$$ 0\to \NS(\bar X)_{\Q_l} \to H^2(\bar X,\Q_l(1))\to V_l\Br(\bar X)\to 0$$
we get
$$%P_2(X,q^{-1}t)=
\det(1-Ft\mid H^2(\bar X,\Q_l(1)))=
\det(1-Ft\mid \NS(\bar X)_{\Q_l})\cdot 
\det(1-Ft\mid V_l\Br(\bar X)).$$

\begin{proposition}\label{braudiv} Assume that $\Br(X)$ is finite and 
let $\sigma'$ be the number of roots $\alpha$ of $L_{tr}(X,t)$ with $v_q(\alpha)=1$.
Then
$$\Br^0(\bar X)\cong (\Q/\Z)^{\sigma'}\oplus (\Q/\Z')^{2g-\sigma'}$$
and
$$ |\Br^0(\bar X)^\Gamma|= L_{tr}(X,q^{-1})q^{\sum_{v(\alpha_i)<1}(1-v(\alpha_i))}.$$
\end{proposition}

\proof %The proof is similar to \cite[Thm. 6.1]{milne}. 
According to the Tate conjecture and \cite[Prop. 5.4, Rem. 5.6]{milne_values}
we have
%$$ \det(1-\varphi t\mid \NS(\bar X)_{\Q_l})=\prod_A(1-\frac{\alpha_i}{q}t);$$
$$\det(1-F t\mid V_l\Br(\bar X))=\prod_{\alpha_i\in B_l}(1-\frac{\alpha_i}{q}t),$$
where %$A$ is the set of $\alpha_j$ which are $q$ times a root of unity, 
$B_l$ is the set of those roots $\alpha_j$ of $P_2(X,t)$ 
which are not $q$ times a root of unity 
if $l\not=p$, and the subset with the further condition 
$v_q(\alpha_i)=1$ if $l=p$. In particular $\dim V_l\Br(\bar X)=2g$
for $l\not=p$ and $\dim V_p\Br(\bar X)=\sigma'$.
The exact sequence
$$ 0\to T_l\Br(\bar X)\to V_l\Br(\bar X)\to \Br^0(\bar X)\{l\}\to 0$$
then proves the statement on $\Br^0(\bar X)$, and the snake Lemma gives
$$
\Br^0(\bar X)^\Gamma = \ker (1-F \mid\Br^0(\bar X))
\cong \bigoplus_l\coker (1-F\mid T_l\Br(\bar X)).$$
%Away from $p$, 
This has order 
\begin{equation}\label{brau}
\prod_l \det(1-F\mid V_l\Br(\bar X))= 
\prod_{l}\prod_{\alpha_i\in B_l}|1-\frac{\alpha_i}{q}|^{-1}_l.
\end{equation}
If $v_q(\alpha_i)=1$, then the product formula shows that 
$\prod_l|1-\frac{\alpha_i}{q}|^{-1}_l=|1-\frac{\alpha_i}{q}|$, but
if $v_q(\alpha_i)\not=1$, then the $p$-part does not appear in 
\eqref{brau}. 
If $v_q(\alpha_i)>1$, then $|1-\frac{\alpha_i}{q}|_p=1$ and we can argue
by the product formula again. If $ v_q(\alpha_i)<1$, then 
$v_q(1-\frac{\alpha_i}{q})=v_q(\alpha_i)-1$, 
hence $|1-\frac{\alpha_i}{q}|^{-1}_p=q^{v_q(\alpha_i)-1}$ and
%, noting that $|x|_p^{-1}=p^{v_p(x)}=q^{v_q(x)}$,
the product formula gives
$$|1-\frac{\alpha_i}{q}|
=q^{v_q(\alpha_i)-1}\prod_{l\not=p}|1-\frac{\alpha_i}{q}|^{-1}_l.$$
The formula in the Proposition follows by recalling that 
the $\alpha_i$ come in complex conjugate pairs with
$v_q(\alpha_i)+v_q(\bar\alpha_i)=2$.
\proofend

\begin{theorem}\label{mainbr}
We have  
$$ |\Br(\bar X)^\Gamma|= 
|\NS(X)_\tor|\cdot L_{tr}(X,q^{-1})\cdot q^{\alpha(X)}.
$$
\end{theorem} 

\proof
By \cite[Remark 7.4]{milne} we have 
$$\alpha(X):=\chi(X,\mathcal O_X)-1+\dim \mathrm{PicVar}_X
=s+\sum_{v(\alpha_i)<1}(1-v(\alpha_i)).$$ 
Moreover, the argument in Proposition \ref{braudiv} shows that 
$H^1(\Gamma, \Br(\bar X)^0)=0$, so that the Theorem follows by 
combining Propositions \ref{brautor} and \ref{braudiv}.
\proofend

\section{Estimating Weil-polynomials}
Consider a Weil-$q^2$-polynomial
$$f(x)=x^{2g}+a_1x^{2g-1}+\cdots +a_1q^{2g-2}x+q^{2g}\in \Z[x]$$
and assume that $\pm q$ is not a root of $f(x)$. Pairing complex conjugate
roots, we can write 
\begin{equation}\label{bequa}
f(x)=\prod_i (x^2+b_ix+q^2),\quad b_i\in \R,  \qquad a_1=\sum_i b_i.
\end{equation}
%and $a_1=-\tr \pi$ if $f(x)$ is irreducible and $\pi$ is a root of $f(x)$.
Let $c_i=2q-b_i$, and consider the polynomial  
$$h(x)=\prod_i (x-c_i).$$
By considering the coefficients of $x^{2g-1},\ldots, x^g$ in $f(t)\in \Z[t]$
one sees that $a_r$ differs from the elementary symmetric function $s_r(b_i)$ 
by integral linear combinations of $s_t(b_i)$ for $t<r$. By induction on $r$
we conclude that $s_r(b_i)\in \Z$. Similarly, the elementary function $s_r(c_i)$
differs from $\pm s_r(b_i)$ by integral linear combinations of the $s_t(b_i)$ 
for $t<r$, hence by induction on $r$ one sees that $h(x)\in \Z[t]$. 
A similar argument shows that  
$f(x)$ is reducible if and only if there is a subset of 
the $b_i$ whose elementary symmetric functions are all integers, if and only if 
$h(x)$ is reducible.

Since $|b_i|<2q$, we have $0<c_i<4q$, hence the roots of $h(x)$ are totally 
positive real algebraic numbers. If $h(x)$ is irreducible, then the roots 
$\xi$ of $h(x)$ have trace 
$$\tr\; \xi =\sum_ic_i %= 2gq+\tr\;\pi
=2qg-a_1.$$
In degrees $g=1, 2$, the irreducible Weil-$q^2$-polynomials with maximal 
$a_1$ are
$$\begin{array}{ll}
f_1(x)=x^2+(2q-1)x+q^2, &h_1(x)=x-1, \\
%and  $a_1> 2qg-q\sqrt e$, and if $g=2$ then 
f_2(x)=x^4+(4q-3)x^3+6q(q-1)x^2+q^2(4q-3)x+q^4, &
h_2(x)=x^2-3x+1.
\end{array}$$
%and $a_1=4q-3>4q-2\sqrt e$.

\begin{remark}
If $f(x)$, hence $h(x)$, is irreducible, then finding the maximal value
for $a_1$
%lower bounds for $\tr\xi$%=2qg-a_i$ 
is equivalent to finding 
the totally positive algebraic integer of degree $g$ with minimal trace. 
This has been extensively discussed in the literature. 

If $g\leq 9$ then $\tr\xi\geq 2g-1$, hence $a_1\leq 2gq-2g+1$, and for 
$g\leq 15$ we have $\tr\xi\geq 2g-2$, hence $a_1\leq 2gq-2g+2$.
Asymptotically, Schur \cite[Satz XI]{schur} showed in 1923 that $\tr\xi>g\cdot \sqrt e$
except for $h_1(x)$ and $h_2(x)$, and Siegel \cite{siegel} showed in 1945 that 
$\tr\xi>g\cdot 1.733$ except for $h_1(x), h_2(x)$, and $h_3(x)=x^3-5x^2+6x+1$.
Currently it is known that $\tr\xi>g\cdot 1.793145$ except for finitely many 
explicitly given exceptions \cite{www}.
\end{remark}

We use Schur's result to give the optimal estimate for $a_1$ valid for all transcendental
Weil-$q^2$ polynomials.

\begin{proposition}\label{c2}
Let $f(x)=x^{2g}+a_1x^{2g-1}+\cdots +a_1q^{2g-2}x+q^{2g}$ be a
Weil-$q^2$-polynomial. %Then $a_1\leq 2gq-2g+1$.
If $f(x)$ does not have $\pm q$ as a root, then $a_1\leq g(2q-1)$.
\end{proposition}

\proof 
Considering each irreducible factor of $f(x)$ separately, 
we can assume that $f(x)$ is irreducible. For $f(x)=f_1(x)$ we have 
$a_1=2q-1$ and for $f(x)=f_2(x)$ we have $a_1=4q-3<g(2q-1)$.
Now $\tr \xi>g\sqrt e $ except for $g_1(x)$ and $g_2(x)$, hence 
$a_1=2gq-\tr\xi < 2gq-g\sqrt e<g(2q-1)$.
\proofend

\begin{theorem}\label{g=2}
If $\deg L_{tr}(X,t)=2g$, then 
$$L_{tr}(X,q^{-1})\leq (4-\frac{1}{q})^g.$$
In particular, for any $X$ with finite Brauer group we  have  
$$|\Br^0(\bar X)^\Gamma|\leq (4q-1)^g.$$
\end{theorem}

\proof Writing $L_{tr}(X,t)=\prod_i(1+b_it+q^2t^2)$  as in \eqref{bequa} and 
comparing geometric and arithmetic means we get
$$ L_{tr}(X,q^{-1})=\prod (2+\frac{b_i}{q}) 
\leq \left(\frac{1}{g}\sum (2+\frac{b_i}{q})\right)^g=(2+\frac{a_1}{gq})^g.$$
By Proposition \ref{c2} we have $a_1\leq g(2q-1)$, so that
$$ L_{tr}(X,q^{-1})\leq (2+\frac{g(2q-1)}{gq})^g= (\frac{4q-1}{q})^g.$$
The final statement follows from this, Proposition \ref{braudiv}, 
and $\sum_{v(\alpha_i)<1}(1-v(\alpha_i))\leq \sum_{v(\alpha_i)<1}1\leq g$.
\proofend

\begin{example}
The estimate in Theorem \ref{g=2} is optimal. 
Assume that the Brauer group of $X$ is finite and that 
$L_{tr}(X,t)=(1+(2q-1)t+t^2q^2)^g$. Since 
the roots of $f_1(t)$ have $v_q=0,2$, we obtain
$\sum_{v(\alpha_i)<1}(1-v(\alpha_i))=g$,
and our calculation gives 
$$\Br^0(\bar X)=(\Q/\Z')^{2g}, \qquad  |\Br^0(\bar X)^\Gamma|=
 (4q-1)^g.$$
See Example \ref{maxab} for an abelian surface with $L_{tr}(X,t)= f_1(t)$.
\end{example}

\section{Abelian surfaces}
We apply our results to prove some results on Brauer groups of abelian surfaces.
Isogeny classes of abelian surfaces over $\F_q$ are classified by their 
Weil-polynomial of degree $4$,  
\begin{equation*}%\label{wpol}
w(x)= x^4-s_1x^3+s_2x^2-qs_1x+q^2
\end{equation*}
where $s_i$ are the elementary symmetric functions in the roots
$\alpha,\beta, \bar\alpha=\frac{q}{\alpha}$, and $\bar\beta=\frac{q}{\beta}$.
From here on we will simply write "abelian variety" instead of "isogeny 
class of abelian varieties".
The reciprocal roots of $P_2(A,T)$ are the $q^2$-Weil numbers
$q,\,q,\,\alpha\beta,\, \alpha\bar \beta,\, \bar\alpha\beta,\, 
\bar\alpha\bar\beta
%q,\,q,\,\alpha\beta,\, \frac{q\alpha}{\beta}=\alpha\bar \beta,\,
%\frac{q\beta}{\alpha}=\bar\alpha\beta,\, 
%\frac{q^2}{\alpha\beta}=\bar\alpha\bar\beta
$, and one can check that this implies that $P_2(A,T)$ equals $(1-qt)^2$ times
\begin{equation}\label{wpol1}
1+(2q-s_2)t+q(2q+s_1^2-2s_2)t^2+q^2(2q-s_2)t^3+q^4t^4.
\end{equation}
We have $\NS(\bar A)_\tor=0$, and Tate's
conjecture is known \cite{tate}, so that
\begin{equation*}%\label{abelformula}
 |\Br(A)|\cdot |\det \NS(A)|=\omega_A |\Br(\bar A)^\Gamma|.
\end{equation*}
We can have $(\rho(\bar A),g)=(6,0), (4,1), (2,2)$, and we describe each case.
See \cite[Cor. 3.12.1, 3.12.2]{suwayui} and \cite[Rem. 3.12.3]{suwayui}
for the classification of Picard numbers. The
data has been assembled with the help of the $L$-function and modular
form data base \cite{page}.

\subsection{\fbox{$\rho(\bar A)=6$}} Then $A$ supersingular, $g=0$, and 
$\Br(\bar A)\cong \bar \F_q$ as Galois modules, so that 
$|\Br(\bar A)^\Gamma|=q$ and 
$$ |\Br(A)|\cdot|\det \NS(A)|=\omega_A q.$$ 
We can have $\tau_A=0,2,4$, hence $\omega_A\in \{1,2,3,4,5,6,8,9,12,16\}$.

\begin{proposition}\label{ss}
For any abelian surface, we have 
$$ |\Br(A)|\cdot|\det \NS(A)|\leq 16 q.$$
This is an equality for the abelian surface $A'$ with Weil-polynomial $(x^2-q)^2$. 
If $\sqrt q\not\in \Z$ then $A'$ is simple, and 
if $\sqrt q\in \Z$, then $A'$ is a product of two non-isogenous 
curves, 
so that $|\Br(A)|=16q$.
\end{proposition}

\proof 
The estimate is clear for supersingular abelian surfaces,
and we will see below that the bound is
smaller for non-supersingular $A$.
For the equality note that $A'$ has $P_2(t)=(1-qt)^2(1+qt)^4$, 
hence $\omega_A=16$. If $\sqrt q\in \Z$, then $A'$ is a product of
the non-isogenous curves corresponding to $(x-\sqrt q)^2$ and $(x+\sqrt q)^2$, 
hence $\det\NS(A)=-1$.
\proofend

The following is a table for possible Weil-polynomials with given $\omega_A$,
using Corollary \ref{exactomega}, \eqref{wpol1}, and the fact that $\rho(A)$ is even. 
The first cases are for all characteristics,
the second group can exist in even degrees extensions of the prime field,
and finally there are some sporadic examples for odd degree
extensions of the prime field in characteristic at most $5$.
$$\begin{array}{r|r|c|c|c}
\omega_A&\tau_A&\frac{P_2(A,t)}{(1-qt)^2}&Weil&for\\
\hline
4&2&(1+qt)^2(1-qt)^2%=1-2q^2t^2+q^4t^4
&x^4+2qx^2+q^2=(x^2+q)^2\\

4&4&(1+qt)^2(1-qt+q^2t^2)%=~1+qt+q^3t^3+q^4t^4
&x^4+qx^2+q^2&\\

8&4&(1+qt)^2(1+q^2t^2)%=1+2qt+2q^2t^2+2q^3t^3+q^4t^4
&x^4+q^2&\\

12&4&(1+qt)^2(1+qt+q^2t^2)%=1+3qt+4q^2t^3+3q^2t^3+q^4t^4 
&x^4-qx^2+q^2\\

16&4&(1+qt)^4%=1+4qt+6q^2t^2+4q^3t^3+q^4t^4
&x^4-2qx^2+q^2=(x^2-q)^2\\
\hline

1&0&(1-qt)^4%=1-4qt+6q^2t^2+4q^3t^3-q^4t^4
&x^4\pm 4\sqrt qx^3+6qx^2\pm 4q\sqrt qx+q^2
&\sqrt q\in\Z\\

1&4&(1-qt+q^2t^2)^2%=
&x^4\pm 3\sqrt qx^3+4qx^2\pm 3q\sqrt qx+q^2&\sqrt q\in\Z\\

1&4&1-q^2t^2+q^4t^4
&x^4\pm \sqrt qx^3+qx^2\pm q\sqrt qx+q^2&\sqrt q\in \Z\\

3&2&(1+qt+q^2t^2)(1-qt)^2%=1-qt-q^3t^3+q^4t^4
&x^4\pm 2\sqrt qx^3+3qx^2\pm 2q\sqrt qx+q^2
&\sqrt q\in \Z\\

4&4&(1+q^2t^2)^2%=1+2q^2t^2+q^4t^4
&x^4\pm 2\sqrt qx^3+2qx^2+2q\sqrt q x+q^2&\sqrt q\in \Z\\

5&4&1+qt+q^2t^2+q^3t^3+q^4t^4&x^4\pm\sqrt qx^3+qx^2\pm q\sqrt qx+q^2&\sqrt q\in \Z\\

9&4&(1+qt+q^2t^2)^2%=1+2qt+3q^2t^3+2q^2t^3+q^4t^4  
&x^4\pm\sqrt qx^3+q\sqrt qx+q^2&\sqrt q\in \Z\\

\hline
2&2&(1+q^2t^2)(1-qt)^2%=1-2qt+2q^2t^2-2q^3t^3-q^4t^4
&x^4\pm 2\sqrt{2q}x^3+4qx^2\pm  2q\sqrt{2q}x^3+q^2 &\sqrt{2q}\in\Z \\

2&4&1+q^4q^4&x^4\pm \sqrt{2q}x^3+2qx^2\pm q\sqrt{2q}x+q^2&\sqrt{2q}\in\Z\\

6&4&(1+qt+q^2t^2)(1+q^2t^2)%=1+qt+2q^2t^2+q^3t^3+q^4t^4
&x^4\pm\sqrt{2q}x^3+qx^2\pm q\sqrt{2q}x+q^2&\sqrt{2q}\in \Z \\

1&2&(1-qt+q^2t^2)(1-qt)^2%=1-3qt+4q^2t^2-3q^3t^3+q^4t^4
&x^4\pm 2\sqrt{3q}x^3+5qx^2\pm 2q\sqrt{3q}x+q^2&\sqrt {3q}\in\Z\\

3&4&(1+qt+q^2t^2)(1-qt+q^2t^2)%=1+q^2t^2+q^4t^4
&x^4\pm \sqrt{3q}x^3+2qx^2\pm q\sqrt{3q}x+q^2& \sqrt{3q}\in\Z\\

1&4&1-qt+q^2t^2-q^3t^3+q^4t^4
&x^4\pm \sqrt{5q}x^3+3qx^2\pm q\sqrt{5q}+q^2&\sqrt{5q}\in\Z \\
\end{array}$$
Note that $P_2(A,t)=(1-qt)^2(1+q^2t^2)(1-qt+q^2t^2)$ does not occur,
and in general not every row occurs for every possible $q$. For example,
for $q=2^2,3^2,5^2,7^2,11^2,13^2$ there are $19$ isogeny classes in the table,
but the actual number is $19,19,15,12,17,9$, respectively. This is due to
the fact that certain Weil-polynomials correspond to higher dimensional abelian
varieties. For example, if $q$ is a square, then 
$x^4+qx^2+q^2=(x^2-\sqrt qx+q)(x^2+\sqrt qx+q)$, which do not correspond to 
elliptic curves for $p\equiv 1\mod 6$.

\begin{example}\label{primef}
Over a prime field $\F_p$, all possible values above appear, and 
we obtain the following decomposition into simple components in the
above table. Here "square" indicates that $A$ is the isogenous to the product 
of the self product of a curve, and "product" indicates that $A$ is isogenous
to the product of two non-isogenous curves. For $p\geq 7$ we get: 
$$\begin{array}{c|c|c|c}
P_1(t)&%&rk\, \NS  & \frac{P_2(A,t)}{(1-pt)^{2}} 
&\omega_A p& |\Br(A)|\\
%&End   \\
\hline
1  +2p t^2 + p^2t^4&square %&4  &(1-pt)^2(1+pt)^2  
&4p&1,4\\
%&M_2(\Q(\sqrt{-p}))\\
1  + p t^2 + p^2t^4&simple %& 2  &(1+pt)^2(1-pt+p^2t^2)
& 4p&1,4\\
%&\Q(\sqrt{-3},\sqrt{p})\\
1 + p^2 t^4        &simple %&2 &(1+pt)^2(1+pt^2)
& 8p&1,4\\
%&\Q(i,\sqrt{2p})\\
1  -p t^2 + p^2t^4   &simple %&2 &(1+pt)^2(1+pt+p^2t^2)
& 12p&1,4 \\
%&\Q(\sqrt{-3},\sqrt{-p})\\
1  -2p t^2 + p^2t^4  &simple %&2  &(1+pt)^4
&  16p &1,4,16\\
%&D/\Q(\sqrt p) \\
\end{array}$$
Over $\F_5$ we have, in addition:
$$\begin{array}{c|c|c|c}
P_1(t)&%&rk\, \NS  & \frac{P_2(A,t)}{(1-pt)^{2}} 
&5\omega_A & |\Br(A)|\\
%&End   \\
\hline
1 \pm 5t+ 15 t^2\pm 25t^3+ 25t^4&simple %& 2&1-5t+25t^2-125t^3+625t^4
&5&1\\
%&\Q(\zeta_5)\\
\end{array}$$
Indeed, for $p\geq 5$ there is only one isogeny class of supersingular elliptic curves 
over $\F_p$, so there is exactly one non-simple abelian surface.
Over $\F_3$, we obtain the following additional 
$4$ isogeny classes and one exceptional product decomposition:
$$\begin{array}{c|c|c|c}
P_1(t)&%&rk\, \NS  & \frac{P_2(A,t)}{(1-3t)^{2}} 
&3\omega_A & |\Br(A)|   \\
\hline
1 \pm 6t+ 15t^2\pm 18t^3+9t^4&square%&4&(1-3t)^2(1-3t+9t^2)
&3&1 \\
1 \pm 3t+ 6 t^2\pm 9t^3+ 9t^4&prod%& 2&(1-3t+9t^2)(1+3t+9t^2)
&9&9 \\
\hline
1  - 3 t^2 + 9t^4 &prod%&2 &(1+3t)^2(1+3t+9t^2)
& 36&36 \\
\end{array}$$
Over $\F_2$, we have the following additional $6$ isogeny classes,
and one exceptional product decomposition:
$$\begin{array}{c|c|c|c}
P_1(t)&%&rk\, \NS & \frac{P_2(A,t)}{(1-2t)^{2}} 
&2\omega_A & |\Br(A)|  \\
\hline
1 \pm 4t+ 8 t^2 \pm 8t^3+ 4t^4&square%& 4 &(1-2t)^2(1+4t^2)
&4 &1,4\\
1 \pm 2t+ 2 t^2 \pm 4t^3+ 4t^4&simple%&2 &(1+4t^2)(1+2t+4t^2)
&12&1,4 \\
1 \pm 2t+ 4 t^2 \pm 4t^3+ 4t^4&prod%&2 &1    +16t^4
& 4 &4\\
\hline
1 + 4t^4           &prod%&2& (1+2t)^2(1+4t^2)
&  16& 16\\
\end{array}$$
\end{example}

A product of two non-isogenous curves has $\NS(A)\cong \Z\times\Z$
with $\det\NS(A)=-1$, and we obtain:

\begin{corollary}\label{ssp}
The largest size of the Brauer group of a supersingular abelian surfaces
over a prime field is $36$, realized by
the abelian surface over $\F_3$ with Weil-polynomial $x^4-3x^2 + 9$.
\end{corollary}

\begin{proposition}
If $q$ is an odd power of the prime $p\not=3$, then for any supersingular 
abelian surface we have 
$$ |\Br(A)|\leq \frac{16q}{p}.$$
If $q$ is a power of $3$, then the $|\Br(A)|\leq 12q$. The bound is sharp
in characteristics $2$ and $3$.
\end{proposition}

\proof 
We know that $|\Br(A)|$ divides $\omega_Aq$, where $\omega_A$ is at most $16$
and has no prime divisor larger than $5$. Moreover, the prime divisor $5$
cannot appear for $q$ an odd power of $p$. 
Taking into account that $|\Br(A)|$ is a square, the result follows for
$p\geq 5$, and for $p=2$ the bound $\frac{16q}{p}$ is the same bound one 
obtains for $\omega_A=8$. 
For $p=3$, the largest value is achieved for $\omega_A=12$. 

Conversely, if $q$ is a power of $3$, then the product of the two non-isogenous
elliptic curves with $L$-polynomials  $1\pm \sqrt{3q}t+qt^2$ has
$L$-polynomial $1-qt^2+q^2t^4$ and $P_2(A,t)=(1-qt)^2(1+qt)^2(1+qt+q^2t^2)$
so that $\omega_A=12$. Since $\det(\NS(A))=1$, we obtain $|\Br(A)|=12q$.

If $q$ is a power of $2$, then the product of the two non-isogenous
elliptic curves with $L$-polynomials  $1\pm \sqrt{2q}t+qt^2$ has
$L$-polynomial $1+q^2t^4$ and $P_2(A,t)=(1-qt)^2(1+qt)^2(1+q^2t^2)$
so that $\omega_A=8$. Since $\det(\NS(A))=1$, we obtain $|\Br(A)|=8q$.
\proofend

If $q\geq 5$, then the abelian surface $A$ with Weil-polynomial
$(x^2-q)^2$ has $P_2(A,t)=(1+qt)^4$ and $\omega_A=16$. 
This abelian surface is simple and the base change of the abelian
surface $A'$ over $\F_p$ with Weil-polynomial $(x^2-p)^2$. For $A'$
we have $|\Br(A')|\cdot |\det(\NS(A'))|=16p$, hence 
$|\det(\NS(A'))|\in\{p,4p,16p\}$.
But the N\'eron-Severi group does not change in the odd degree extension 
$\F_q/\F_p$. However, we do not know how to determine the $2$-part.

\subsection{\fbox{$\rho(\bar A)=4$}} 
Then $L_{tr}(A,t)$ has degree $2$ and $\bar A$ is the self-product of a
(necessarily ordinary) elliptic curve, 
and $A$ has $p$-rank $2$. The slopes of Frobenius of $P_2(A,t)$ 
are $0,1,1,1,1,2$ and 
$\Br(\bar A)\cong (\Q/\Z')^2$. 
The Frobenius eigenvalues with slope one are $q$ times a root of unity,
and $|\Br(\bar A)^\Gamma|\leq 4q-1$ by Theorem \ref{g=2}.

\subsubsection{$\rho(A)=2, \tau_A=2, \omega_A\in \{1,2,3,4\}$.} Then $A$ is simple
or the product of two non-isogenous elliptic curves (which become isogenous
after a finite base extension).

\begin{proposition}\label{largest}
If $\rho(A)=2$ and $\rho(\bar A)=4$, then 
$$|\Br(A)|\cdot |\det \NS(A)|\leq 16q-4.$$
If $\sqrt q\in \Z$ then there is an abelian surface with 
$$|\Br(A)|=4(2\sqrt q-1)^2,$$ 
and if $\sqrt q\not\in \Z$ and $p\not|\lfloor 2\sqrt q\rfloor$, 
then there is an abelian surface with 
$$|\Br(A)|=4\lfloor 2\sqrt q\rfloor^2.$$ 
These are the largest possible order for a Brauer group.
\end{proposition}

\proof
By Theorem \ref{g=2} we obtain the upper bound
$$ |\Br(A)|\cdot |\det \NS(A)|=\omega_A\cdot |\Br(\bar A)^\Gamma|\leq 4(4q-1)$$
and $\omega_A$ divides the product on the left hand side.
%We have $\sum_{v_q(\alpha_i)<1}(1-v_q(\alpha_i))=1$,
%hence $ |\Br(A)|\cdot |\det \NS(A)|=8q+s_1^2-4s_2$ by \eqref{rho2}.
If $\rho(A)=2$, then substituting $q^{-1}$ into \eqref{wpol1} we obtain
using the original Artin-Tate formula that
\begin{equation}\label{rho2}
|\Br(A)|\cdot |\det \NS(A)|=8q+s_1^2-4s_2.
\end{equation}
If $\sqrt q\in \Z$ consider the polynomial 
$$x^4-(2q-4\sqrt q+1)x^2+ q^2=
(x^2-(2\sqrt q-1) x+q)(x^2+(2\sqrt q-1) x+q).$$
This is the product of
two non-isogenous ordinary elliptic curves, hence $|\det \NS(A)|=1$.
Substituting $s_1$ and $s_2$ into the above we obtain 
$ |\Br(A)|=4(2\sqrt q-1)^2$. 

If $\sqrt q\not\in \Z$ and consider the polynomial 
$$x^4+(2q-\lfloor 2\sqrt q\rfloor^2)x^2+ q^2=
(x^2-\lfloor 2\sqrt q\rfloor x+q)(x^2+\lfloor 2\sqrt q\rfloor x+q).$$
Since we assume that $p\not| \lfloor 2\sqrt q\rfloor$, this is  
the product of two non-isogenous ordinary elliptic curves, 
hence $|\det \NS(A)|=1$.
Substituting $s_1$ and $s_2$ into \eqref{rho2} we obtain 
$ |\Br(A)|=4\lfloor 2\sqrt q\rfloor^2$. Note that since $\sqrt q\not \in\Z$, 
$\lfloor 2\sqrt q\rfloor^2\leq 4q-1$, in accordance to our estimate. 

To show that this is the largest order of a Brauer group, 
assume that $A$ has a larger Brauer group. If
$|\det \NS(A)|\not=1$ or $\omega_A\leq 2$, then $|\Br(A)|\leq 2(4q-1)=8q-2$
and for all $q>2$ we have 
$$4\lfloor 2\sqrt q\rfloor^2>4(2\sqrt q-1)^2= 16q-16\sqrt q+4>8q-2$$
and for $q=2$ we still have $4\lfloor 2\sqrt q\rfloor^2=16>14=8q-2$.
If $|\det \NS(A)|=1$ and $\omega_A=3$, then the order of the Brauer
group is divisible by $9$ and $|\Br(A)|\leq 3(4q-1)=12q-3$. 
As above one checks that for $q>11$, 
$4(2\sqrt q-1)^2>12q-3$. In the remaining cases $q\leq 11$, the value in
the proposition is at least as large as the largest square no exceeding
$12q-3$ which is divisible by $9$. For example, for $q=7$, the value in the 
proposition is $100>12q-3=81$. 
%But $(2\lfloor 2\sqrt q\rfloor)^2>(2(2\sqrt q-1))^2= 16q-16\sqrt q+4$, 
%and this is larger than $12q-3$ if $q\geq 13$, and for $q\leq 11$ one easily 
%checks explicitly that the values in the proposition are the largest squares
%below $12q-3$. 
Thus $|\det \NS(A)|=1$ and $\omega_A=4$, so that 
$|\Br(A)|$ is divisible by $4$. Now it suffices to observe
that $4(2\sqrt q-1)^2$ and $4\lfloor 2\sqrt q\rfloor^2$ are the largest squares
below $4(4q-1)$, respectively.
\proofend

%If $A$ is the product of two non-isogenous curves, then $\det \NS(A)=1$, and 
%if $s_1=0$, then $P_2(X,t)=(1-2t)^2(1+2t)^2\cdot L_{tr}(X,t)$,
%hence $\omega=4$. 

\begin{example}\label{maxab}
Let  $A$ be the abelian surface over $\F_q$ with Weil-polynomial 
$x^4-(2q-1)x^2+q^2$. Since $16q-4$ is never a square, $A$ is simple by 
\cite[Thm. 2.9]{nart}, and after a degree $2$ extension, $A_{\F_{q^2}}$ is the 
product of two isogenous elliptic curves with Weil polynomial 
$x^2-(2q-1)x+q^2$. We have 
$$P_2(A,t)=(1+(2q-1)t+q^2t^2)(1+qt)^2(1-qt)^2,$$ 
hence $\omega_A=4$ and
$$ |\Br(A)|\cdot |\det \NS(A)|=|\Br(\bar A)^\Gamma|\cdot \omega_A =16q-4.$$ 
\end{example}

\subsubsection{$\rho(A)=4, \tau_A=0$.} Then $A$ itself
is isogenous to a self-product of an ordinary elliptic curve $E$ over $\F_q$, 
$\omega_A=1$,  
$$|\coker\alpha|=|\det \NS(A)|=|\disc(\End(E))|$$
and 
$$|\Br(A)|\cdot |\det \NS(A)|=|\Br(\bar A)^\Gamma|\leq  4q-1. $$

\begin{example}
Let  $A$ be the abelian surface over $\F_q$ with Weil-polynomial
$(x^2\pm x+q)^2$. Then  
$P_2(A,t)=(1+(2q-1)t+q^2t^2)(1-qt)^4$ 
hence $$|\Br(A)|\cdot |\det \NS(A)|=4q-1,$$
the maximum possible value.
Now $x^2\pm x+q$ has roots $\pi=\frac{-1\pm\sqrt{-(4q-1)}}{2}$, 
hence $\End(E)$ is the maximal order with discriminant $-(4q-1)$. 
We obtain 
%Since $|\det \NS(A)|$ is a multiple of this, we obtain 
$\Br(A)=0$ and $|\det \NS(A)|=4q-1$.
\end{example}

\begin{remark}
In \cite[Prop. 5.6]{yui} it is claimed that $\Br(A)=0$ for a product
of ordinary elliptic curves, but this is not true because the colimit of 
$\Br(A\times_{\F_q}\F_{q^r})$ is isomorphic to $\Br(\bar A)\cong (\Q/\Z')^2$.
\end{remark}

\subsection{\fbox{$\rho(\bar A)=2$}}  Then $L_{tr}(A,t)$ has degree $4$, 
$\tau_A=0$ and $\omega_A=1$. Either $\bar A$ is simple or the 
product of two non-isogenous curves.

\begin{proposition}\label{largest2}
If $\rho(\bar A)=2$, then 
$$|\Br(A)|\cdot |\det \NS(A)|\leq (4\sqrt q-1)^2.$$
There is an abelian surface with 
$$|\Br(A)|=(2\lfloor 2\sqrt q\rfloor -1)^2,$$ 
and this is the largest possible order for the Brauer group of an abelian
surface with $\rho(\bar A)=2$.
\end{proposition}

\proof
As in Proposition \ref{largest} we have 
%$\sum_{v_q(\alpha_i)<1}(1-v_q(\alpha_i))=1$, hence 
$ |\Br(A)|\cdot |\det \NS(A)|=8q+s_1^2-4s_2$. % by \eqref{rho2}.
By \cite[Lemma 2.1]{nart}
we have %$|s_1|\leq 4\sqrt q$ and 
$-4s_2\leq 8q-8|s_1|\sqrt q$, hence 
$$|\Br(A)|\cdot |\det \NS(A)|\leq 16q-8|s_1|\sqrt q+s_1^2=(4\sqrt q-|s_1|)^2.$$
If $s_1=0$, then the Weil-polynomial of $A$ is biquadratic and 
$A$ is the self-product of an elliptic curve after a degree $2$ extension,
contradicting $g=2$. Thus $(4\sqrt q-1)^2$ is an upper bound. 
Conversely consider the polynomial 
%$$x^4\pm x^3-2(q-\sqrt q)x^2\pm qx+q^2=(x+\sqrt q)^2(x^2-(2\sqrt q-1)x+q).$$
$$x^4+ x^3+(2q-\lfloor 2\sqrt q\rfloor^2+\lfloor 2\sqrt q\rfloor)x^2+ qx+q^2=
(x^2+\lfloor 2\sqrt q\rfloor x+q)(x^2-(\lfloor 2\sqrt q\rfloor -1)x+q).$$
This is the Weil-$q$ polynomial of the product of two non-isogenous 
elliptic curves, at least one of which is ordinary, hence $|\det \NS(A)|=1$.
Substituting $s_1$ and $s_2$ into the above we obtain 
$ |\Br(A)|=(2\lfloor 2\sqrt q\rfloor -1)^2$.
This finished the proof if $\sqrt q\in \Z$. For $\sqrt q\not\in \Z$ and $|s_1|=1$,
note that
$8q+1-4s_2$ is odd, so that the next possible order of the Brauer group is 
$(2\lfloor 2\sqrt q\rfloor +1)^2$. But $\lfloor 2\sqrt q\rfloor>2\sqrt q-1$
so that $2\lfloor 2\sqrt q\rfloor+ 1>4\sqrt q-1$, a contradiction.
On the other hand, if $|s_1|>1$ then the next larger order of the Brauer 
group is $(2\lfloor 2\sqrt q\rfloor)^2>(4\sqrt q-2)^2$, a contradiction as well.
\proofend

\iffalse
The situation for $q$ an odd prime power is more complicated because to
have $g=2$ we need either $p$-rank one, or equivalently 
$s_2$ divisible by $p$,
or $A$ has to be a product of non-isonegenous ordinary curves. 
%(this is automatic in the optimal value above for $q$ a square).

\begin{remark}
If $q$ is an odd power of $p$, then there is an abelian surface with 
$|\Br(A)|=(2\lfloor 2\sqrt{q}\rfloor-1)^2$.
\end{remark}
\fi

\subsubsection{$p$-rank one} Any abelian surface of $p$-rank one
has $\rho(\bar A)=2$. The Frobenius eigenvalues of $P_2(A,t)$ are 
$\frac{1}{2}, \frac{1}{2},1,1, \frac{3}{2}, \frac{3}{2}$,
hence $\Br(\bar A)\cong (\Q/\Z')^4$.
%$$ |\Br(A)|\cdot |\det \NS(A)|=|\Br(\bar A)^\Gamma|=qL_{tr}(A,q^{-1})\leq 16q-8.$$

\iffalse
\begin{example}(2.2b-a) Let  $A$ be the simple abelian variety over 
$\F_2$ with Weil-polynomial $x^4+x^3+2x+4$. Then 
$P_2(A,t)=  (1+4x+ 10x^2+16x^3+16x^4)(1-2t)$ and 
$$ |\Br(A)|\cdot |\det \NS(A)|=|\Br(\bar A)^\Gamma|=17.$$
Since $|\Br(A)|$ is a square we get $\Br(A)=\{0\}$ and 
$|\det \NS(A)|=17$. Going through all $8$ isogeny classes 
of abelian surfaces of $p$-rank $1$， one obtains $1,9,17$, so that $17$ is
the maximum.
\end{example}
\fi

\begin{example}
If $p>3$, then there is only one supersingular elliptic curves over $\F_p$ with Weil-polynomial
$x^2+p$, hence a non-simple abelian surface is the product of this curve and the ordinary elliptic
curve with Weil-polynomial $x^2+ax+p, 0<|a|\leq 2\sqrt{p}$. Since $\det \NS(X)=1$ for the product
of non-isogenous curves, we obtain that for every $0<a\leq  2\sqrt{p}$
there is an abelian surface of $p$-rank one with $|\Br(X)|=a^2$. The situation for simple abelian 
surfaces is more complicated.
\end{example}

\subsubsection{$p$-rank two} 
An abelian variety of $p$-rank two is in this category if $\bar A$ simple 
or the product of two non-isogenous ordinary curves.
The Frobenius eigenvalues of $P_2(A,t)$ are $0,1,1,1,1,2$,  
hence $\Br(\bar A)\cong (\Q/\Z)^2\oplus (\Q/\Z')^2$.

\iffalse
\begin{example}
$(2.9.ab-al)$ Let  $A$ be the absolutely simple abelian surface over $\F_9$
with L-polynomial $1 - x - 11 x^{2} - 9 x^{3} + 81 x^{4}$. 
One can calculate that $P_2(A,t)=(1-7t+9\cdot 41t^2-9^2\cdot 7t^3+9^4t^4)(1-9t)^2$, 
hence 
$$ |\Br(A)|\cdot |\det \NS(A)|=|\Br(\bar A)^\Gamma| =117=3^2\cdot 13.$$
%Indeed, the endomorphism algebra has $\Q(\sqrt{13})$ as its unique
%subfield of degree $2=\rho(A)$, so that $|\det \NS(A)|=13$.
\end{example}
\fi

\end{document}